\documentclass[11pt,a4paper]{article}

\usepackage[T1]{fontenc}
\usepackage[latin1]{inputenc}

\usepackage{amsmath}
\usepackage{amssymb}
\usepackage{amsthm}
\usepackage{amsfonts}
\usepackage{mathrsfs}

\newcommand{\Ann}[0]{\operatorname{Ann}}

\newcommand{\cohdim}[0]{\operatorname{cd}}
\newcommand{\Hom}[0]{\operatorname{Hom}}
\newcommand{\Ass}[0]{\operatorname{Ass}}
\newcommand{\Coass}[0]{\operatorname{Coass}}

\newcommand{\Supp}[0]{\operatorname{Supp}}

\newcommand{\Spec}[0]{\operatorname{Spec}}

\newcommand{\height}[0]{\operatorname{height}}
\newcommand{\cd}[0]{\operatorname{cd}}

\newtheorem{theorem}{Theorem}[section]

\newtheorem{corollary}[theorem]{Corollary}

\newtheorem{remark}[theorem]{Remark}

\newtheorem{question}[theorem]{Question}

\title{On the associated primes of Matlis duals of local cohomology modules II}

\author{Michael Hellus}

\begin{document}

\maketitle

\begin{abstract}

In continuation of \cite{hellus05} we study associated primes of Matlis duals of local cohomology modules (MDLCM). We combine ideas from Helmut Z\"oschinger on coassociated primes of arbitrary modules with results from \cite{hellus05}, \cite{hellus08}, \cite{hellus_stueckrad08}, \cite{hellus_stueckrad_to_appear} and obtain partial answers to questions which were left open in \cite{hellus05}. These partial answers give further support for conjecture $(*)$ from \cite{hellus05} on the set of associated primes of MDLCMs. In addition, and also inspired by ideas from Z\"oschinger, we prove some non-finiteness results of local cohomology.

\end{abstract}

\section{Introduction}

Let $I$ be an ideal of a local, noetherian ring $R$. By $H^l_I$ we denote the $l$-th local cohomology functor supported on $I$, by $E$ a fixed $R$-injective hull of the residue field of $R$ and by $D$ the Matlis dual functor $D:=\Hom_R(\_ ,E)$ from $(R-mod)$ to $(R-mod)$.

Suppose that one has $H^l_I(R)=0$ for $l\neq c$ ($c$ is necessarily the height of $I$ then). Assume that a regular sequence $x_1,\ldots ,x_c$ in $I$ is given. It was shown in the author's Habilitationsschrift (\cite[Cor. 1.1.4]{hellus06}) that $I$ is a set-theoretic complete intersection defined by the $x_i$ if and only if the $x_i$ form a $D(H^c_I(R))$-(quasi)regular sequence. This gives strong motivation to study the associated primes of $D(H^c_I(R))$. It is this study which we started in \cite{hellus05} and which we continue here.

The simplest case is $R=k[[X_1,\ldots ,X_n]]$ and $I=(X_1,\ldots ,X_c)R$, where $k$ is a field, the $X_i$ are indeterminates and $0\leq c\leq n$. The case $c=n$ is easy; the case $c=n-1$ is non-trivial and was completely solved in \cite[Theorem 2.5]{hellus_stueckrad08}, see also \cite{hellus05}. The next case is $c=n-2$, where the following is known (\cite[Theorem 2.2.1]{hellus05} and \cite[Theorem 1.3(ii),(v)]{hellus_stueckrad08}):

\begin{itemize}

\item \[ p\in \Ass_R(D(H^{n-2}_{(X_1,\ldots ,X_{n-2})}(R)))\Rightarrow \height p\in \{ 0,1,2\} .\]

\item \[ \{ 0\} \in \Ass_R(D(H^{n-2}_{(X_1,\ldots ,X_{n-2})}(R))).\]

\item If $\height p=2$:
\[ p\in \Ass_R(D(H^{n-2}_{(X_1,\ldots ,X_{n-2})}(R)))\iff \sqrt{p+(X_1,\ldots ,X_{n-2})}=\sqrt{(X_1,\ldots X_n)}.\]

\item If $\height p=1$: $P$ is generated by a prime element $p$ of $R$: $P=pR$. If $p\not\in (X_1,\ldots ,X_{n-2})$, then
\[ pR\in \Ass_R(D(H^{n-2}_{(X_1,\ldots ,X_{n-2})}(R))).\]

\end{itemize}

It is natural to ask next

\begin{question}

\label{qu_ht_one}Which height-one prime ideals, i.~e. which $P=pR$, where $p$ is an (arbitrary) prime element of $R$, are in $\Ass_R(\underbrace{D(H^{n-2}_{(X_1,\ldots ,X_n-2)}(R))}_{=:D})$?

\end{question}

This question is open (but note that some {\it very special} height one prime ideals in $\Ass_R(D)$ where found in \cite[Cor. 4.3.1]{hellus06}). The main goal of this paper is to show that in many cases the answer to question \ref{qu_ht_one} is {\it positive}; in particular, it is positive if $k$ is countable and $p$ is a polynomial contained in $(X_{n-1},X_n)R$. In fact, our two main results, theorem \ref{allg_Fall} and theorem \ref{Polynomfall}, are both a little more general, see section \ref{Associated prime ideals} for the precise statements. An example which is by no means trivial and where question \ref{qu_ht_one} has a positive answer is given by $p=X_{n-1}X_1+X_nX_2$ (if $n\geq 4$, of course). This example follows from theorem \ref{Polynomfall}.

The results in section \ref{Associated prime ideals} give some indication that conjecture $(*)$ from \cite[section 1]{hellus05} (which says in this situation that \[ \Ass_R(D)=\{ p|H^{n-2}_{(X_1,\ldots ,X_{n-2})}(R/p)\neq 0\} \]) holds, because in the situation of theorem \ref{allg_Fall} one has $H^{n-2}_I(R/(a,b)R)\neq 0$ and, a fortiori, $H^{n-2}_I(R/pR)\neq 0$; in this context, see also \cite[Theorem 1.1]{hellus05}.

In section \ref{non-fin_prop} we prove some non-finiteness properties of local cohomology modules: It is very well-known that top local cohomology modules are almost never finitely generated, see e.~g. \cite[Remark 2.5]{hellus08_2} for a quick proof using the Nakayama lemma. In fact a stronger statement holds: No quotient of a top local cohomology module is finite (corollary \ref{no_finite_quotient}), and we do not even have to assume that the module $M$ whose local cohomology we consider must be finite. One even has that top local cohomology modules have no coatomic quotients (theorem \ref{no_koatomar_quotient}; a module is {\it coatomic} if every proper submodule is contained in a maximal one).

Helmut Z\"oschinger's work on coatomic modules and coassociated prime ideals (e.~g., \cite{zoeschinger80}, \cite{zoeschinger86}, \cite{zoeschinger88}) is essential for both sections of this paper.

\section{Associated prime ideals}

\label{Associated prime ideals}

By 'countable' we shall mean either finite or 'infinite countable'.

\begin{theorem}

\label{allg_Fall}Let $k$ be a countable field, $R$ a domain and a local $k$-algebra essentially of finite type, $n:=\dim R\geq 4$, $I\subseteq R$ an ideal, $\height I=n-2=\cd I$. Assume that there exist $a,b\in R$ such that $(a,b)R$ is prime and $a,b$ define a system of parameters for $R/I$, let $p\in (a,b)R$ be a prime element. Then \[pR\in \Ass_RD(H^{n-2}_I(R)).\]

\end{theorem}

{\it Proof. }Obviously $R$ has only countably many prime ideals (as any algebra of finite type over $k$ has only countably many (prime) ideals). By \cite[Theorem 2.1]{hellus_stueckrad_to_appear} there exist infinitely many prime ideals $q$ which contain $p$ and which are associated to $D(H^{n-2}_I(R))$. For each such $q$ one has in particular $0\neq \Hom_R(R/q,D(H^{n-2}_I(R)))\buildrel (*_1)\over =D(H^{n-2}_I(R)\otimes_RR/q)\buildrel (*_2)\over =D(H^{n-2}_I(R/q))$ ($(*_1)$: Hom-Tensor adjointness, $(*_2)$: Right exactness of $H^{n-2}_I$) and hence $\height(q)\leq 2$. As, therefore, all these $q$ have either height one (in which case $q$ equals $pR$) or height two, their intersection is $pR$ (the height of this intersection is necessarily one, as infinitely many pairwise different $q$s are intersected). It follows that the intersection of all associated prime ideals of $\Hom_R(R/pR,D(H^{n-2}_I(R)))$ is $pR$. By \cite[Lemma 3.1]{zoeschinger86}, the associated prime ideals of $D(H^{n-2}_I(R))$ are precisely the coassociated prime ideals of $H^{n-2}_I(R)$. \cite[Folgerung 1.5 and Lemma 3.1]{zoeschinger88} imply that each prime ideal minimal over $pR$ is associated to $H^{n-2}_I(R)$. But $pR$ is prime and hence we get $pR\in \Coass_R(H^{n-2}_I(R))=\Ass_RD(H^{n-2}_I(R))$.\hfill $\square $

\begin{theorem}

\label{Polynomfall}Let $k$ be a field, $X_1,\ldots ,X_n$ indeterminates, $n\geq 4$. Set $R=k[[X_1,\ldots ,X_n]]$ and $I=(X_1,\ldots ,X_{n-2})R$. Let $p\in (X_{n-1},X_n)R$ be a prime element that has $pR\cap R_0\neq 0$, where $R_0:=k_0[X_1,\ldots ,X_n]_{(X_1,\ldots ,X_n)}$ and where $k_0$ is a countable subfield of $k$ (e.~g. the prime subfield of $k$). Then
\[ pR\in \Ass_R(D(H^{n-2}_I(R))).\]

\end{theorem}

{\it Proof.} $pR\cap R_0$ has height at most one, by our hypothesis it must hence have the form $p_0R_0$ for some prime element $p_0\in R_0$ (note that prime elements are non-zero by definition). As $k_0$ is countable, we get from theorem \ref{allg_Fall}
\[ p_0R_0\in \Ass_{R_0}(D(H^{n-2}_{(X_1,\ldots ,X_{n-2})R_0}(R_0)))\]
(here $D$ is taken with respect to $R_0$, of course). By \cite[Lemma 3.1]{zoeschinger86}, $p_0R_0\in \Coass_{R_0}(\underbrace{H^{n-2}_{(X_1,\ldots ,X_{n-2})R_0}(R_0)}_{=:H})$. That means there exists an Artinian quotient $H\twoheadrightarrow H/B$ of $H$ that has
\[ \Ann_{R_0}(H/B)=p_0R_0.\]
The $R$-module
\[ (H/B)\otimes_{R_0}R\buildrel{R/R_0\text{ faithfully flat}}\over=(H\otimes_{R_0}R)/(B\otimes_{R_0}R)\]
is a quotient of $H\otimes_{R_0}R$ and is Artinian (because its support is zero-dimensional and its socle
\[ \Hom_R(R/m,(H/B)\otimes_{R_0}R)=\Hom_{R_0}(R_0/(X_1,\ldots ,X_n),H/B)\otimes_{R_0}R\] has finite vector space-dimension); furthermore, by faithful flatness of $R/R_0$, its annihilator is
\[ \Ann_R((H/B)\otimes_{R_0}R)=p_0R.\]
By Matlis duality, $D((H/B)\otimes_{R_0}R)$ is a finitely generated $R$-submodule of $D(H^{n-2}_I(R))$ with annihilator
\[ \Ann_R(D((H/B)\otimes_{R_0}R))=\Ann_R((H/B)\otimes_{R_0}R)=p_0R.\]
The prime ideal $pR$ is minimal over $p_0R$, therefore we get
\[ pR\in \Ass_R(D((H/B)\otimes_{R_0}R))\subseteq D(H^{n-2}_I(R)).\]\hfill $\square $
\begin{remark}

\begin{itemize}

\item In the situation of theorem \ref{Polynomfall} one can quickly show that $\{ 0\} \in \Ass_RD(H^{n-2}_I(R))$ using the following arguments (this case was already known, with a different proof, see \cite[Lemma 2.1.1]{hellus05}): The intersection of all coassociated prime ideals of $H^{n-2}_I(R)$ equals the radical of $\Ann_RH^{n-2}_I(R)$ (this follows from \cite[Satz 1.2 and Folgerung 1.3]{zoeschinger88}, because $0=H^{n-2}_I(R/(X_1,\ldots ,X_n))=H^{n-2}_I(R)\otimes_RR/(X_1,\ldots ,X_n)R$, i.~e. one has $(X_1,\ldots ,X_n)H^{n-2}_I(R)=H^{n-2}_I(R)$); but the endomorphism ring of $H^{n-2}_I(R)$ is $R$, by \cite[Theorem 2.2 (iii)]{hellus08}; in particular, $\Ann_RH^{n-2}_I(R)=0$. Therefore, using the argument from the proof of theorem \ref{allg_Fall}, one concludes $\{ 0\} \in \Ass_R(D(H^{n-2}_I(R)))$.

\end{itemize}

\end{remark}

It seems natural to ask

\begin{question}

In the situation of theorem \ref{Polynomfall}, is it true that
\[ pR\in \Ass_R(D(H^{n-2}_I(R)))\]
holds for {\rm every} prime element $p\in (X_{n-1},X_n)R$?

\end{question}

\begin{question}

Does conjecture $(*)$ hold in this context, i.~e. is it true that

\[ \Ass_R(D(H^{n-2}_I(R)))=\{ p\in \Spec R|H^{n-2}_I(R/p)\neq 0\} ?\]

\end{question}

With respect to prime ideals of height two or zero both questions have positive answer, this was explained in the introduction. The results in this paper say that both questions have at least {\it often} a positive answer for height one prime ideals.

\section{Non-finiteness properties}

\label{non-fin_prop}

Whenever, over a local, complete ring $(R,m)$, a given local cohomology module $H$ has infinitely many coassociated prime ideals (this is often the case: \cite[Theorem 3.1.3 (ii), (iii)]{hellus06}), $H$ is neither finitely generated (because if it was, then $D(H)$ would be Artinian and hence one would have $\Ass_RDH=\{ m\} $) nor Artinian (because if it was then $\Ass_R(D(H))$ would be finite). This trivial remark is generalized.

\begin{remark}

\label{basic_facts_coatomic_m}Over the noetherian ring $R$, the coatomic modules are closed under taking quotients, submodules and extensions, see \cite[section 1]{zoeschinger80}. It is clear that every finitely generated $R$-module is coatomic and that every coatomic, Artinian module has finite length. Furthermore, localizations of coatomic modules are coatomic (over the localized ring), see \cite[section 1, Folgerung 2]{zoeschinger80}.

\end{remark}

\begin{theorem}

\label{no_koatomar_quotient}Let $R$ be a noetherian ring, $M$
an $R$-module and $I$ an ideal of $R$ such that $1\leq c:=\cohdim (I,M)=\cohdim(I,R/\Ann_R(M))<\infty $ (without further assumption one would have only $\cohdim (I,M)\leq \cohdim(I,R/\Ann_R(M))$ in general). Then the top local cohomology module $H^c_I(M)$ has no non-zero coatomic quotient.

\end{theorem}

{\it Proof. }If $H^c_I(M)$ had a non-zero, coatomic quotient $H^c_I(M)/U$, then, by localizing in an arbitrary $p\in \Supp_R(H^c_I(M)/U)$, we would get a non-zero, coatomic (remark \ref{basic_facts_coatomic_m}) quotient of $H^c_I(M)_p=H^c_{IR_p}(M_p)$. Therefore, we may replace $R$ by $R_p$ and assume that $(R,m)$ is local (note also that one has $c=\cohdim (I,M)=\cohdim(IR_p,M_p)$).

Assume to the contrary that $H/U$ is a non-zero, coatomic quotient
of $H:=H^c_I(M)$ for some submodule $U$ of $H$. In particular there exists a maximal submodule $U'$ of $H$ containing $U$. Being a simple module, $H/U'$ is ismomorphic to $R/m$.

On the other hand, $D(H/U')$ is naturally a submodule of $D(H)$ and it is also isomorphic to $R/m$. But $m$ is not associated to $D(H/U')\subseteq D(H)$ (because otherwise
\[ 0\neq
\Hom_R(R/m,D(H))=D(H^c_I(M)\otimes_R(R/m))\buildrel (\dag)\over=D(H^c_I(M/mM))=0,\]
contradiction; for $(\dag)$ one works over the ring $R/\Ann_R(M)$ and uses the fact that $H^c_{I(R/\Ann_R(M))}$ is right exact on $R/\Ann_R(M)$-modules). Therefore, no such quotient $H/U$ exists and the
theorem is proven.\par\hfill $\square $

Note that in the formulation of theorem \ref{no_koatomar_quotient} (as well as in the subsequent corollary \ref{no_finite_quotient})
it is not required that $M$ is finitely generated.

\begin{corollary}

\label{no_finite_quotient}Let $I$ be an ideal of a noetherian ring
and let $M$ be an $R$-module such that $1\leq c:=\cohdim (I,M)=\cohdim(I,R/\Ann_R(M))<\infty $. Then
$H^c_I(M)$ has no non-zero finitely generated quotient.

\end{corollary}

\begin{remark}

The proof of the preceding theorem actually shows that in the given situation the top local cohomology module is {\rm radikalvoll} (see e.~g. \cite{zoeschinger80} for this terminology: By definition, a module is {\rm radikalvoll} if it has no maximal submodule).

\end{remark}

As an application of theorem \ref{no_koatomar_quotient} we get immediately an improvement of \cite[Cor. 1.1.4]{hellus06} (recall that a sequence $(x_1,\ldots ,x_n)$ in a local ring $R$ is {\it filter regular} on the $R$-module $M$ if, for each $i$, the kernel of the multiplication map $M/(x_1,\ldots ,x_{i-1})M\buildrel x_i\over \to M/(x_1,\ldots ,x_{i-1})M$ is Artinian) see e.~g. \cite{schenzel78} and \cite{stueckrad86}):

\begin{theorem}

Let $(R,m)$ be a noetherian, local ring, $I$ a proper ideal of $R$,
$h\in \mathbb N$ and $\underline f=f_1,\dots ,f_h\in I$ an
$R$-regular sequence. The following statements are equivalent:

\begin{enumerate}

\item $\sqrt {\underline fR}=\sqrt I$.

\item \label{quasi}$H^l_I(R)=0$ for every $l>h$ and the sequence $\underline
f$ is quasi-regular on $D(H^h_I(R))$.

\item $H^l_I(R)=0$ for every $l>h$ and the sequence $\underline
f$ is regular on $D(H^h_I(R))$.

\item \label{filter} $H^l_I(R)=0$ for every $l>h$ and the sequence $\underline
f$ is filter regular on $D(H^h_I(R))$.

\end{enumerate}

\end{theorem}

{\it Proof. }Because of \cite[Cor. 1.1.4]{hellus06} it suffices to show that \ref{filter}. implies \ref{quasi}: Assume that $h\geq 1$, $H^l_I(R)=0$ for
every $l>h$ and that $\underline f=f_1,\ldots ,f_h\in I$ is a
filter regular sequence on $D(H^h_I(R))$. In particular, the kernel $K$
of the multiplication map
\[ D(H^h_I(R))\buildrel f_1\over \to D(H^h_I(R))\]
is Artinian. But
$K=\Hom_R(R/f_1,D(H^h_I(R)))=D(H^h_I(R)\otimes_R(R/f_1R))$ and hence
the quotient module $H^h_I(R)\otimes_R(R/f_1)$ is a finitely
generated $\hat R$-module. It follows from theorem
\ref{no_koatomar_quotient} that $K=0$. But then we have
$D(H^{h-1}_I(R/f_1R))=D(H^h_I(R))\otimes _R(R/f_1)R)$ by an easy
argument with exact sequences. Now it is clear that the claim
follows by induction on $h$.\hfill $\square $

Michael Hellus, Universit\"at Leipzig, Fakult\"at für Mathematik und Informatik, PF 10 09 20, 04009 Leipzig, Germany

E-Mail: hellus@math.uni-leipzig.de

\end{document}